\documentclass[12pt,A4paper]{article}
\usepackage{epsfig}
\usepackage{graphicx}
\usepackage{graphpap}
\usepackage{amsmath}
\usepackage{amsfonts}
\usepackage{amssymb}
\usepackage{placeins}
\def\baselinestretch{2.0}
\def\newcommand\abstractname {abstract}

\def\abstractname {}

\begin{document}

\baselineskip=24pt
\begin{center}
\large{THE ASYMPTOTIC EFFICIENCY OF IMPROVED PREDICTION INTERVALS}
\end{center}

\bigskip
\begin{center}
{\large {\large B}{\normalsize Y} {\large P}{\normalsize AUL} {\large K}{\normalsize ABAILA$^*$ AND}
{\large K}{\normalsize HRESHNA} {\large S}{\normalsize YUHADA}}
\end{center}

\begin{center}
\large{\textit{La Trobe University and Institut Teknologi Bandung}}
\end{center}

\def\baselinestretch{1.0}

\vspace{11cm}

\noindent * Author to whom correspondence should be addressed.

\noindent Department of Mathematics and Statistics, La Trobe University, Victoria 3086, Australia.

\noindent e-mail: P.Kabaila@latrobe.edu.au

\noindent Facsimile: 3 9479 2466

\noindent Telephone: 3 9479 2594

\newpage

\baselineskip=21pt

\textbf{Abstract.} 
Barndorff-Nielsen and Cox (1994, p.319) modify an estimative prediction
limit to obtain an improved prediction limit with better coverage 
properties. 
Kabaila and Syuhada (2008) present a simulation-based approximation to this
improved prediction limit, which avoids the extensive algebraic manipulations
required for this modification. 
We present a modification of an estimative prediction interval, analogous to the 
Barndorff-Nielsen and Cox modification,
to obtain an improved
prediction interval with better coverage properties. We also present an analogue, for the 
prediction interval context,
of this simulation-based approximation. 
The parameter
estimator on which the estimative and improved prediction limits and intervals are based is assumed
to have the same asymptotic distribution as the (conditional) maximum
likelihood estimator. 
The improved prediction limit and interval depend on the asymptotic 
conditional bias of this estimator.
This bias can be very sensitive to very small changes in the estimator.
It may require considerable effort to find this bias.
We show, however, that the improved
prediction limit and interval have asymptotic efficiencies that are functionally independent
of this bias. Thus, 
improved prediction limits and intervals obtained using the Barndorff-Nielsen and Cox 
type of methodology can conveniently 
be based on the (conditional) maximum likelihood estimator, whose asymptotic conditional
bias is given by the formula of Vidoni (2004, p.144). Also, 
improved prediction limits and intervals obtained using Kabaila and Syuhada type
approximations have asymptotic efficiencies that are independent of the estimator on which
these intervals are based. 

\bigskip


\def\baselinestretch{2.0}

\newpage

\baselineskip=24pt

\begin{center}
{\bf 1. INTRODUCTION}
\end{center}

\noindent Suppose that $\{Y_t\}$ is a discrete-time stochastic
process with probability distribution determined by the parameter
vector $\theta$, 
where the $Y_t$ are continuous random variables. Also suppose that $\{Y^{(t)}\}$ is a Markov
process, where $Y^{(t)} = \big (Y_{t-p+1},\ldots,Y_t \big)$. For
example, $\{Y_t\}$ may be an AR($p$) process or an ARCH($p$)
process. The available data is $Y_1, \ldots, Y_n$. Suppose that we
are concerned with $k$-step-ahead prediction where $k$ is a
specified positive integer. 
Also suppose that $\widehat \Theta$ is an estimator of $\theta$ with the same
asymptotic distribution as the (conditional) maximum likelihood estimator.
We note that there are many possible choices for $\widehat \Theta$. For example,
for a stationary Gaussian AR(1) model, commonly-used estimators of the 
autoregressive parameter include least-squares,
Yule-Walker and Burg estimators. 
We use lower case to denote observed values of random vectors. For example,
$y^{(n)}$ denotes the observed value of the random vector $Y^{(n)}$. 
We also use the Einstein summation notation that repeated indices are implicitly
summed over.

Firstly, suppose that our aim is to find an upper prediction
limit $z(Y_1,\ldots, Y_n)$, for $Y_{n+k}$, 
such that it has coverage probability conditional
on $Y^{(n)}=y^{(n)}$ equal to $1-\alpha$ i.e. such that
\begin{equation*}
P_\theta \big ( Y_{n+k} \le z(Y_1,\ldots, Y_n) \, \big| \,
Y^{(n)} = y^{(n)} \big ) = 1 - \alpha
\end{equation*}
for all $\theta$ and $y^{(n)}$. The desirability of a prediction limit or
interval having coverage probability $1-\alpha$ conditional
on $Y^{(n)}=y^{(n)}$ has been noted by a number of authors. In the context of an
AR($p$) process, this has been noted by Phillips (1979), Stine (1987), Thombs and
Schucany (1990), Kabaila (1993), McCullough (1994), He (2000), Kabaila and He (2004)
and Vidoni (2004). In the context of an ARCH($p$) process, this has been noted by
Christoffersen (1998), Kabaila (1999), Vidoni (2004) and Kabaila and Syuhada (2008).

Define $z_{\alpha}(\theta, y^{(n)})$ by the requirement that 
$P_\theta \big ( Y_{n+k} \le z_{\alpha}(\theta, y^{(n)}) \, \big| \,
Y^{(n)} = y^{(n)} \big ) = 1 - \alpha$ for all $\theta$ and $y^{(n)}$.
The estimative $1-\alpha$ prediction limit is defined to be $z_{\alpha}(\widehat \Theta, Y^{(n)})$.
This prediction limit may not have adequate coverage probability properties
unless $n$ is very large. It may be shown that the coverage probability of 
$z_{\alpha}(\widehat \Theta, Y^{(n)})$ conditional on $Y^{(n)}=y^{(n)}$ differs from $1-\alpha$ 
by $O(n^{-1})$. 
Barndorff-Nielsen and Cox (1994, p. 319) 
modify (using a procedure clarified by Vidoni (2004)
and described in detail by Kabaila and Syuhada, 2008, Section 2)
the estimative prediction limit to obtain an improved
prediction limit with better coverage properties. This improved
limit, denoted $z_{\alpha}^+(\widehat \Theta, Y^{(n)})$ and described in Section 3 
of the present paper, has
coverage probability
conditional on $Y^{(n)} = y^{(n)}$ that differs from $1 - \alpha$ by $O(n^{-3/2})$.
The algebraic manipulations needed to obtain this improved prediction limit
are feasible only for the simplest time series models. To avoid these 
manipulations, Kabaila and Syuhada (2008) propose a simulation-based approximation to this
improved prediction limit.

In the present paper, 
we extend these results to prediction intervals as follows.
In Section 4, we show that the estimative $1-\alpha$ prediction interval has
coverage probability
conditional on $Y^{(n)} = y^{(n)}$ that differs from $1 - \alpha$ by $O(n^{-1})$.
In this section, we also present a modification of an estimative 
$1-\alpha$ prediction interval, analogous to the 
Barndorff-Nielsen and Cox (1994, p.319) modification of an estimative prediction limit,
to obtain an improved $1-\alpha$ prediction interval with better coverage properties. 
We show that this improved $1-\alpha$ prediction interval has coverage probability
conditional on $Y^{(n)} = y^{(n)}$ that differs from $1 - \alpha$ by $O(n^{-3/2})$.
To avoid the extensive algebraic manipulations required to find this improved 
prediction interval, we propose a simulation-based approximation to this
interval, analogous to Kabaila and Syuhada (2008)
simulation-based approximation to the improved prediction limit.

The improved $1-\alpha$ prediction limit and interval are obtained 
from the estimative $1-\alpha$ prediction limit and interval, respectively, 
using a correction that includes the asymptotic bias of $\widehat \Theta$ conditional on 
$Y^{(n)} = y^{(n)}$. Kabaila and Syuhada (2007, Section 4) present an example showing that
this bias can be very sensitive to small changes in the estimator $\widehat \Theta$.
A further illustration of this fact is provided in Section 2 of the present paper.
It may require considerable effort to find the asymptotic bias of $\widehat \Theta$ conditional on 
$Y^{(n)} = y^{(n)}$. In Sections 3 and 4 we show, however, that the improved $1-\alpha$ 
prediction limit and interval have asymptotic efficiencies that are functionally independent of this
bias. This has the following two consequences. Firstly, if the improved prediction limit or interval is 
obtained algebraically using the Barndorff-Nielsen and Cox (1994, p. 319) methodology or its analogue, 
respectively, then the estimative and improved prediction limits and intervals can conveniently
be based on the (conditional) maximum likelihood estimator. This is because the 
asymptotic conditional bias of this estimator can be found using the very convenient formula of 
Vidoni (2004, p.144). Secondly, improved prediction limits and intervals obtained using
a Kabaila and Syuhada (2008) type simulation-based approximation have asymptotic efficiency 
that is independent of the estimator $\widehat \Theta$, on which both the estimative and improved
$1-\alpha$ prediction limits and intervals are based. Note that we assume throughout this paper 
that $\widehat \Theta$ has the same asymptotic distribution as the 
(conditional) maximum likelihood estimator.

\bigskip
\begin{center}
2. SENSITIVITY OF THE ASYMPTOTIC CONDITIONAL BIAS TO
SOME SMALL CHANGES IN THE ESTIMATOR
\end{center}

\noindent Consider a stationary zero-mean Gaussian AR(1) process $\{Y_t\}$
satisfying
$ Y_t = \rho \, Y_{t-1} + \varepsilon_t,$
for all integer $t$, where $|\rho|<1$ and the $\varepsilon_t$ are
independent and identically $N(0,\sigma^2)$ distributed. Note that
$\varepsilon_t$ and $(Y_{t-1},Y_{t-2},\ldots)$ are independent for
each $t$. The available data is $Y_1,Y_2,\ldots,Y_n$.

The least-squares estimator
$
\widehat{\rho} = \sum^{n}_{t=2} \, Y_t \,
Y_{t-1}/\sum^{n-1}_{t=1} \, Y_t^2  
$
is obtained by maximizing the log likelihood function
conditional on $Y_1=y_1$. The Yule-Walker estimator
$
\widehat{\rho}_{YW} = \sum^{n}_{t=2} \, Y_t \,
Y_{t-1}/\sum_{t=1}^n \, Y_t^2 
$
differs by a very small amount from $\widehat \rho$. However, as proved by
Shaman and Stine (1988), $E(\widehat \rho - \rho) = -2 \rho n^{-1}+ \cdots$
and $E(\widehat \rho_{YW} - \rho) = -3 \rho n^{-1}+ \cdots$. This illustrates the 
great sensitivity of the asymptotic (unconditional) bias of an estimator of $\rho$
to some small changes in this estimator.

We illustrate the great sensitivity of the asymptotic bias 
of an estimator of $\rho$
conditional
on the last observation to some small changes in this estimator as
follows. Define the estimator
$
\widetilde{\rho}=\sum^n_{t=2} \, Y_t
Y_{t-1}/\sum^n_{t=2} \, Y^2_t, 
$
which is obtained by maximizing the log likelihood function
conditional on $Y_n =y _n$. This log likelihood function is found
using the backward representation of the process:
$
Y_t = \rho \, Y_{t+1} + \eta_t,
$
for all integer $t$, where the $\eta_t$ are independent and
identically $N(0,\sigma^2)$ distributed. Note that $\eta_t$ and
$(Y_{t+1}, Y_{t+2}, \ldots )$ are independent.

The estimators $ \widehat{\rho}$ and $\widetilde{\rho}$ differ by only a small amount. 
They have the same asymptotic (unconditional) bias, since
$E(\widehat \rho - \rho) = -2 \rho n^{-1}+ \cdots$ and
$E(\widetilde \rho - \rho) = -2 \rho n^{-1}+ \cdots$.
Yet their
asymptotic biases conditional on $Y_n=y_n$ are quite different.
These asymptotic conditional biases are described as follows.
\begin{align*}
E \big( \widehat{\rho}-\rho|Y_n=y_n \big) &= \big( y_n^2 \, (1-\rho^2)
\rho \, (\sigma^2)^{-1} - 3 \, \rho \big) \, n^{-1} + \cdots \\
E \big( \widetilde{\rho}-\rho|Y_n=y_n \big) &= -2 \, \rho \, n^{-1} +
\cdots
\end{align*}
These expressions for asymptotic bias may be obtained using the
formula for the asymptotic conditional bias of the maximum
likelihood estimator described by Vidoni
(2004, p. 144).

\newpage

\begin{center}
{3. EFFICIENCY RESULT FOR IMPROVED PREDICTION LIMITS}
\end{center}

\noindent Let $F(\, \cdot \, ; \theta, y^{(n)})$ denote the cumulative distribution function
of $Y_{n+k}$, conditional on $Y^{(n)} = y^{(n)}$. Also, let 
$f(\, \cdot \, ; \theta, y^{(n)})$ denote the probability density function corresponding to
this cumulative distribution function. Assume, as do Barndorff-Nielsen and Cox (1994) and Vidoni (2004),
that 
\begin{align}
\label{bias}
  E_{\theta} \big( \widehat \Theta - \theta \, | \,  Y^{(n)} = y^{(n)} \big)
  &=  b(\theta, y^{(n)})n^{-1} + \cdots \\
\label{inverse_exp_info_matrix}
  E_{\theta} \big((\widehat \Theta - \theta) (\hat{\Theta} - \theta)^T
    \, \big | \,  Y^{(n)} = y^{(n)} \big)
  &= i^{-1}(\theta) + \cdots
\end{align}
where $i(\theta)$ denotes the expected information matrix. 
We assume that every element of $i(\theta)$ is $O(n^{-1})$.

Define $H_{\alpha}(\theta|y^{(n)}) = P_{\theta} \big (Y_{n+k} \le z_{\alpha}(\widehat \Theta,y^{(n)}) 
\, \big | \, Y^{(n)} = y^{(n)} \big)$, which is the conditional coverage probability of the 
$1-\alpha$ estimative prediction limit. Using the fact that the distribution of $Y_{n+k}$ given
$(Y_1, \ldots, Y_n)=(y_1, \ldots, y_n)$ depends only on $y^{(n)}$, it may be shown that
$H_{\alpha}(\theta|y^{(n)}) = E_{\theta} \big (F(z_{\alpha}(\widehat \Theta,y^{(n)}); \theta,y^{(n)})
\, \big | \, Y^{(n)} = y^{(n)} \big)$. Now define
$G_{\alpha} (\widehat \Theta; \theta|y^{(n)}) = F(z_{\alpha}(\widehat \Theta,y^{(n)}); \theta,y^{(n)})$.
Thus $H_{\alpha}(\theta|y^{(n)}) = E_{\theta} \big(G_{\alpha} (\widehat \Theta; \theta|y^{(n)}) \, \big | \, 
Y^{(n)} = y^{(n)} \big)$. We now use the stochastic expansion
\begin{align}
\notag
\label{G_expansion}
G_{\alpha} (\widehat \Theta; \theta|y^{(n)}) = &G_{\alpha} (\theta; \theta|y^{(n)}) 
+ \frac{\partial G_{\alpha} (\widehat \theta; \theta|y^{(n)})}{\partial \widehat \theta_i} \Bigg|_{\hat{\theta} = \theta} 
(\widehat \Theta_i - \theta_i)\\
&+ \frac{1}{2} 
\frac{\partial^2 G_{\alpha}(\widehat \theta; \theta|y^{(n)})}{\partial \widehat \theta_r \partial \widehat \theta_s}
\Bigg|_{\hat{\theta}=\theta} (\widehat \Theta_r - \theta_r) (\widehat \Theta_s - \theta_s)
+ \cdots
\end{align}
By the definition of $z_{\alpha}(\theta,y^{(n)})$, $G_{\alpha} (\theta; \theta|y^{(n)})=1-\alpha$. Thus
$H_{\alpha}(\theta|y^{(n)}) = 1-\alpha + c_{\alpha}(\theta,y^{(n)}) n^{-1} + \cdots$ where
\begin{equation}
\label{c}
c_{\alpha}(\theta,y^{(n)}) n^{-1} =   n^{-1} \, \frac{\partial G_{\alpha}(\widehat \theta;\theta|y^{(n)})}{\partial \widehat \theta_i} 
\Bigg|_{\hat{\theta} = \theta} \,
b(\theta, y^{(n)})_i  
+ \frac{1}{2} 
\frac{\partial^2 G_{\alpha}(\widehat \theta; \theta|y^{(n)})}{\partial \widehat \theta_r \partial \widehat \theta_s}
\Bigg|_{\hat{\theta} = \theta} i^{rs}  
\end{equation}
where $b(\theta, y^{(n)})_i$ denotes the $i$th component of the vector $b(\theta, y^{(n)})$ and 
$i^{rs}$ denotes the $(r,s)$th element of the inverse of the expected information matrix $i(\theta)$.
In other words, the conditional coverage probability of the estimative $1-\alpha$ upper prediction limit
$z_{\alpha}(\widehat \Theta,Y^{(n)})$ is $1-\alpha + O(n^{-1})$.

Define
\begin{equation}
\label{d}
d_{\alpha}(\theta,y^{(n)}) 
= - \frac{c_{\alpha}(\theta,y^{(n)}) n^{-1}}{f(z_{\alpha}(\theta,y^{(n)});\theta,y^{(n)})}.
\end{equation}
The improved $1-\alpha$ prediction limit described by Barndorff-Nielsen and Cox (1994, p.319) is
\begin{equation*}
z_{\alpha}^+ (\widehat \Theta,Y^{(n)}) = z_{\alpha} (\widehat \Theta,Y^{(n)}) + d_{\alpha}(\widehat \Theta,Y^{(n)}).
\end{equation*}
The conditional coverage probability of this improved prediction limit is 
$P_{\theta} \big (Y_{n+k} \le z_{\alpha}^+ (\widehat \Theta,y^{(n)}) 
\, \big | \, Y^{(n)} = y^{(n)} \big) =
E_{\theta} \big ( F(z_{\alpha}^+ (\widehat \Theta,y^{(n)});\theta,y^{(n)}) \, \big | \, Y^{(n)} = y^{(n)} \big)$.
We now use the expansion
\begin{align*}
F \big(z_{\alpha}^+ (\widehat \Theta,y^{(n)});\theta,y^{(n)} \big) &= F\big(z_{\alpha} (\widehat \Theta,y^{(n)});\theta,y^{(n)}\big)
+ f \big(z_{\alpha}(\widehat \Theta,y^{(n)});\theta,y^{(n)}\big) \, d_{\alpha}(\widehat \Theta,y^{(n)}) + \cdots \\
&=G_{\alpha} (\widehat \Theta;\theta|y^{(n)})
+ f \big(z_{\alpha}(\theta,y^{(n)});\theta,y^{(n)} \big) \, d_{\alpha}(\theta,y^{(n)}) + \cdots
\end{align*}
Thus 
\begin{align*}
&P_{\theta}(Y_{n+k} \le z_{\alpha}^+ (\widehat \Theta,Y^{(n)}) \, | \, 
Y^{(n)} = y^{(n)})  \\
&=
H_{\alpha}(\theta|y^{(n)}) + f \big(z_{\alpha}(\theta,y^{(n)});\theta,y^{(n)} \big) \, d_{\alpha}(\theta,y^{(n)}) + \cdots \\
&= 1 - \alpha + O(n^{-3/2})
\end{align*}
Note that the improved prediction limit $z_{\alpha}^+ (\widehat \Theta,Y^{(n)})$ may be found algebraically using 
\eqref{c} and \eqref{d}. When these algebraic manipulations become too complicated, 
the method of Kabaila and Syuhada (2008) may be used. For any given $\theta$, these authors estimate 
$P_{\theta} \big (Y_{n+k} \le z_{\alpha}(\widehat \Theta,y^{(n)})\, | \, 
Y^{(n)} = y^{(n)}) - (1-\alpha)$
by Monte Carlo simulation and use this estimate as an approximation to $c_{\alpha}(\theta,y^{(n)}) n^{-1}$
(which appears in \eqref{d}).
In Kabaila and Syuhada (2008), the formula for $r(\omega, y^{(n)})$ should be
$n^{-1} c(\omega, y^{(n)})/f(z(y^{(n)},\omega);\omega| y^{(n)})$ instead of 
$c(\omega, y^{(n)})/f(z(y^{(n)},\omega);\omega| y^{(n)})$, so that  
$d(\omega, y^{(n)})= n^{-1} c(\omega, y^{(n)})$, to order $n^{-1}$.

We measure the asymptotic efficiency of the improved prediction limit \newline
$z_{\alpha}^+ (\widehat \Theta,Y^{(n)})$ by examining the asymptotic expansion of 
$E_{\theta} (z_{\alpha}^+ (\widehat \Theta,Y^{(n)}) \, | \, Y^{(n)} = y^{(n)})$. In other
words, this asymptotic efficiency is a function of $\theta$ and $y^{(n)}$. 
This measure of asymptotic efficiency is consistent with the general guidelines put forward
by Kabaila and Syuhada (2007) for comparing the efficiencies of prediction intervals.
Using $G_{\alpha} (\widehat \Theta; \theta|y^{(n)}) = F(z_{\alpha}(\widehat \Theta,y^{(n)}); \theta,y^{(n)})$, we find that
\begin{align*}
d_{\alpha}(\theta,y^{(n)}) &= -  n^{-1} \, \frac{\partial z_{\alpha}(\theta,y^{(n)})}{\partial \theta_i} \,
b(\theta, y^{(n)})_i  \\
&- \left ( \frac{f^{\prime}(z_{\alpha}(\theta,y^{(n)});\theta,y^{(n)})}{2 f(z_{\alpha}(\theta,y^{(n)});\theta,y^{(n)})}
\frac{\partial z_{\alpha}(\theta,y^{(n)})}{\partial \theta_r} \frac{\partial z_{\alpha}(\theta,y^{(n)})}{\partial \theta_s}
+ \frac{1}{2} \frac{\partial^2 z_{\alpha}(\theta,y^{(n)})}{\partial \theta_r \partial \theta_s} \right) i^{rs}  
\end{align*}
Now, $z_{\alpha}^+ (\widehat \Theta,y^{(n)})$
is equal to
\begin{align*}
&z_{\alpha} (\theta,y^{(n)}) + 
\frac{\partial z_{\alpha} (\theta,y^{(n)})}{\partial \theta_i} (\widehat \Theta_i - \theta_i)
+ \frac{1}{2} \frac{\partial^2 z_{\alpha} (\theta,y^{(n)})}{\partial \theta_r \partial \theta_s}
(\widehat \Theta_r - \theta_r)(\widehat \Theta_s - \theta_s) \\
&+ d_{\alpha}(\theta,y^{(n)}) + \cdots 
\end{align*}
Thus $E_{\theta} \big(z_{\alpha}^+ (\widehat \Theta,y^{(n)})\, \big| \, Y^{(n)} = y^{(n)} \big)$ is equal to
\begin{equation*}
z_{\alpha} (\theta,y^{(n)}) 
- 
\frac{f^{\prime}(z_{\alpha} (\theta,y^{(n)}); \theta,y^{(n)})}
{2 f(z_{\alpha} (\theta,y^{(n)}); \theta,y^{(n)})} \,
\frac{\partial z_{\alpha} (\theta,y^{(n)})}{\partial \theta_r} \,
\frac{\partial z_{\alpha} (\theta,y^{(n)})}{\partial \theta_s} \, i^{rs} + \cdots
\end{equation*}
We see that the asymptotic conditional bias $b(\theta, y^{(n)})n^{-1}$ does not enter into this 
expression. This has the following two consequences. Firstly, if the improved prediction limit is 
found algebraically using \eqref{c} and \eqref{d} then we can use that estimator $\widehat \Theta$ whose 
asymptotic conditional bias is easiest to find. Usually, this will be the (conditional) maximum likelihood estimator
whose asymptotic conditional bias can be found using the formula of
 Vidoni (2004, p.144). Secondly, if
the simulation-based method of Kabaila and Syuhada (2008) is used 
then we know that the asymptotic efficiency of 
the improved $1-\alpha$ prediction limit is 
independent of the estimator $\widehat \Theta$, on which the estimative $1-\alpha$ prediction limit
is based. Note that we assume throughout this paper that $\widehat \Theta$ has the same asymptotic distribution as the 
(conditional) maximum likelihood estimator.

\bigskip
\begin{center}
{4. RESULTS FOR IMPROVED PREDICTION INTERVALS}
\end{center}

\noindent Suppose that our aim is to find a prediction interval 
$\big [ \ell(Y_1,\ldots,Y_n), \, u(Y_1,\ldots,Y_n) \big ]$ for $Y_{n+k}$, such that it has
coverage probability conditional on $Y^{(n)}=y^{(n)}$ equal to $1-\alpha$ i.e. such that
\begin{equation*}
P_\theta \big ( Y_{n+k} \in \big [ \ell(Y_1,\ldots,Y_n), \, u(Y_1,\ldots,Y_n) \big ] \, \big| \,
Y^{(n)} = y^{(n)} \big ) = 1 - \alpha
\end{equation*}
for all $\theta$ and $y^{(n)}$. As in Section 3, define 
$F(\, \cdot \, ; \theta, y^{(n)})$ and $f(\, \cdot \, ; \theta, y^{(n)})$
to be the cumulative distribution function and probability density function (respectively)
of $Y_{n+k}$, conditional on $Y^{(n)} = y^{(n)}$. Suppose that $f(\, \cdot \, ; \theta, y^{(n)})$
is a continuous unimodal function for all $y^{(n)}$ and $\theta$.

Define $\ell_{\alpha}(\theta,y^{(n)})$ and $u_{\alpha}(\theta,y^{(n)})$ by the requirements that
$f(\ell_{\alpha}(\theta,y^{(n)}) ; \theta, y^{(n)}) =
f(u_{\alpha}(\theta,y^{(n)}) ; \theta, y^{(n)})$ and
\begin{equation*}
P_\theta \big ( Y_{n+k} \in \big [ \ell_{\alpha}(\theta,y^{(n)}), \, u_{\alpha}(\theta,y^{(n)}) \big ] \, \big| \,
Y^{(n)} = y^{(n)} \big ) = 1 - \alpha
\end{equation*}
for all $\theta$ and $y^{(n)}$. If $\theta$ is known then 
$\big [ \ell_{\alpha}(\theta,y^{(n)}), \, u_{\alpha}(\theta,y^{(n)}) \big ]$ 
is the shortest prediction interval for $Y_{n+k}$, 
having coverage probability $1-\alpha$ 
conditional on $Y^{(n)}=y^{(n)}$.
We define the estimative $1-\alpha$ prediction interval to be 
\begin{equation*}
I_{\alpha}(\widehat \Theta, Y^{(n)}) = 
\big [ \ell_{\alpha}(\widehat \Theta,Y^{(n)}), \, u_{\alpha}(\widehat \Theta,Y^{(n)}) \big ].
\end{equation*}
Assume that \eqref{bias} and \eqref{inverse_exp_info_matrix} hold true.

Define $H_{\alpha}(\theta|y^{(n)}) = P_{\theta} \big (Y_{n+k} \in I_{\alpha}(\widehat \Theta,y^{(n)}) 
\, \big | \, Y^{(n)} = y^{(n)} \big)$, which is the conditional coverage probability of the 
$1-\alpha$ estimative prediction interval. Using the fact that the distribution of $Y_{n+k}$ given
$(Y_1, \ldots, Y_n)=(y_1, \ldots, y_n)$ depends only on $y^{(n)}$, it may be shown that
$H_{\alpha}(\theta|y^{(n)}) = E_{\theta} \big(G_{\alpha} (\widehat \Theta; \theta|y^{(n)}) \, \big | \, 
Y^{(n)} = y^{(n)} \big)$, where we define 
$G_{\alpha} (\widehat \Theta; \theta|y^{(n)}) = 
F(u_{\alpha}(\widehat \Theta,y^{(n)}); \theta,y^{(n)}) -
F(\ell_{\alpha}(\widehat \Theta,y^{(n)}); \theta,y^{(n)})$.
We now use the expansion \eqref{G_expansion}. By the definition of 
$\ell_{\alpha}(\theta,y^{(n)})$ and $u_{\alpha}(\theta,y^{(n)})$,
$G_{\alpha} (\theta; \theta|y^{(n)})=1-\alpha$. Thus
$H_{\alpha}(\theta|y^{(n)}) = 1-\alpha + c_{\alpha}(\theta,y^{(n)}) n^{-1} + \cdots$ where
$c_{\alpha}(\theta,y^{(n)}) n^{-1}$ is given by \eqref{c}.
In other words, the conditional coverage probability of the estimative $1-\alpha$ upper prediction interval
$I_{\alpha}(\widehat \Theta, Y^{(n)})$ is $1-\alpha + O(n^{-1})$.

Suppose that $d_{\alpha}^{\ell}(\theta,y^{(n)})$ and $d_{\alpha}^u(\theta,y^{(n)})$ are both
$O(n^{-1})$ for every $\theta$ and $y^{(n)}$. Also suppose that 
\begin{equation}
\label{sum_2_d}
d_{\alpha}^{\ell}(\theta,y^{(n)}) + d_{\alpha}^u(\theta,y^{(n)})
= - \frac{c_{\alpha}(\theta,y^{(n)}) n^{-1}}{f(u_{\alpha}(\theta,y^{(n)});\theta,y^{(n)})}.
\end{equation}
Note that we could replace $f(u_{\alpha}(\theta,y^{(n)});\theta,y^{(n)})$ 
in the denominator of the expression on the
right-hand side by 
$f(\ell_{\alpha}(\theta,y^{(n)});\theta,y^{(n)})$, since
$f(\ell_{\alpha}(\theta,y^{(n)});\theta,y^{(n)}) =
f(u_{\alpha}(\theta,y^{(n)});\theta,y^{(n)})$.
The improved $1-\alpha$ prediction interval is
\begin{equation*}
I_{\alpha}^+(\widehat \Theta, Y^{(n)})  = 
\big [ \ell_{\alpha}(\widehat \Theta,Y^{(n)}) - d_{\alpha}^{\ell}(\widehat \Theta,Y^{(n)}), 
\, u_{\alpha}(\widehat \Theta,Y^{(n)}) + d_{\alpha}^u(\widehat \Theta,Y^{(n)}) \big ].
\end{equation*}
The conditional coverage probability of this improved prediction interval is 
\begin{align*}
&P_{\theta} \big (Y_{n+k} \in I_{\alpha}^+(\widehat \Theta, y^{(n)}) 
\, \big | \, Y^{(n)} = y^{(n)} \big) \\
&=
E_{\theta} \big ( F(u_{\alpha} (\widehat \Theta,y^{(n)})+ d_{\alpha}^u(\widehat \Theta,Y^{(n)});\theta,y^{(n)}) \\
& \phantom{123456789012}
- F(\ell_{\alpha} (\widehat \Theta,y^{(n)}) - d_{\alpha}^{\ell}(\widehat \Theta,Y^{(n)});\theta,y^{(n)}) 
\, \big | \, Y^{(n)} = y^{(n)} \big)
\end{align*}
We now use the stochastic expansion
\begin{align*}
&F(u_{\alpha} (\widehat \Theta,y^{(n)})+ d_{\alpha}^u(\widehat \Theta,Y^{(n)});\theta,y^{(n)}) 
- F(\ell_{\alpha} (\widehat \Theta,y^{(n)}) - d_{\alpha}^{\ell}(\widehat \Theta,Y^{(n)});\theta,y^{(n)}) \\
&= F\big(u_{\alpha} (\widehat \Theta,y^{(n)});\theta,y^{(n)}\big)
+ f \big(u_{\alpha}(\widehat \Theta,y^{(n)});\theta,y^{(n)}\big) \, d_{\alpha}^u (\widehat \Theta,y^{(n)}) \\
&\phantom{123}- F\big(\ell_{\alpha} (\widehat \Theta,y^{(n)});\theta,y^{(n)}\big)
+ f \big(\ell_{\alpha}(\widehat \Theta,y^{(n)});\theta,y^{(n)}\big) \, d_{\alpha}^{\ell} (\widehat \Theta,y^{(n)}) + \cdots \\
&=G_{\alpha} (\widehat \Theta;\theta|y^{(n)})
+ f \big(u_{\alpha}(\theta,y^{(n)});\theta,y^{(n)} \big) \, 
(d_{\alpha}^{\ell}(\theta,y^{(n)})+d_{\alpha}^u(\theta,y^{(n)})) \,  + \cdots
\end{align*}
Thus 
\begin{align*}
&P_{\theta}(Y_{n+k} \in I_{\alpha}^+(\widehat \Theta, y^{(n)}) \, | \, 
Y^{(n)} = y^{(n)})  \\
&=
H_{\alpha}(\theta|y^{(n)}) + f \big(z_{\alpha}(\theta,y^{(n)});\theta,y^{(n)} \big) \, 
\big(d_{\alpha}^{\ell}(\theta,y^{(n)})+d_{\alpha}^u(\theta,y^{(n)}) \big) + \cdots \\
&= 1 - \alpha + O(n^{-3/2})
\end{align*}
Note that the improved prediction interval $I_{\alpha}^+(\widehat \Theta, Y^{(n)})$ may be 
found algebraically using \eqref{c} and \eqref{sum_2_d}.
When these algebraic manipulations become too complicated, 
a simulation-based method, similar to that described by Kabaila and Syuhada (2008)
for prediction intervals, may be used. For any given $\theta$, we estimate
$P_{\theta}(Y_{n+k} \in I_{\alpha}(\widehat \Theta, y^{(n)})\, | \, 
Y^{(n)} = y^{(n)}) - (1-\alpha)$  by Monte Carlo simulation and use this estimate 
as an approximation to $c_{\alpha}(\theta,y^{(n)}) n^{-1}$
(which appears in \eqref{sum_2_d}).

We measure the asymptotic efficiency of the improved prediction interval
$I_{\alpha}^+(\widehat \Theta, Y^{(n)})$ by examining the asymptotic expansion of 
$E_{\theta} \big (\text{length of } I_{\alpha}^+(\widehat \Theta, y^{(n)}) \, \big| \, $ $ Y^{(n)} = y^{(n)} \big)$. In other
words, this asymptotic efficiency is a function of $\theta$ and $y^{(n)}$. 
Using $G_{\alpha} (\widehat \Theta; \theta|y^{(n)}) = 
F(u_{\alpha}(\widehat \Theta,y^{(n)}); \theta,y^{(n)}) -
F(\ell_{\alpha}(\widehat \Theta,y^{(n)}); \theta,y^{(n)})$, we find that
\begin{equation*}
\frac{\partial G_{\alpha} (\widehat \theta; \theta|y^{(n)})}{\partial \widehat \theta_i} \Bigg|_{\hat{\theta} = \theta} 
= f \big(u_{\alpha}(\theta,y^{(n)});\theta,y^{(n)} \big) 
\left ( \frac{\partial u_{\alpha}(\theta,y^{(n)})}{\partial \theta_i}
-\frac{\partial \ell_{\alpha}(\theta,y^{(n)})}{\partial \theta_i} \right)
\end{equation*}
and 
\begin{align*}
\frac{\partial^2 G_{\alpha}(\widehat \theta; \theta|y^{(n)})}{\partial \widehat \theta_r \partial \widehat \theta_s}
\Bigg|_{\hat{\theta}=\theta} 
= &f^{\prime} \big(u_{\alpha}(\theta,y^{(n)});\theta,y^{(n)} \big) 
\frac{\partial u_{\alpha}(\theta,y^{(n)})}{\partial \theta_r}
\frac{\partial u_{\alpha}(\theta,y^{(n)})}{\partial \theta_s} \\
&+ f \big(u_{\alpha}(\theta,y^{(n)});\theta,y^{(n)} \big) 
\frac{\partial^2 u_{\alpha}(\theta,y^{(n)})}{\partial \theta_r \partial \theta_s} \\
&-f^{\prime} \big(\ell_{\alpha}(\theta,y^{(n)});\theta,y^{(n)} \big) 
\frac{\partial \ell_{\alpha}(\theta,y^{(n)})}{\partial \theta_r}
\frac{\partial \ell_{\alpha}(\theta,y^{(n)})}{\partial \theta_s} \\
&- f \big(\ell_{\alpha}(\theta,y^{(n)});\theta,y^{(n)} \big) 
\frac{\partial^2 \ell_{\alpha}(\theta,y^{(n)})}{\partial \theta_r \partial \theta_s}
\end{align*}
We now substitute these expressions into \eqref{c} and \eqref{sum_2_d}, to obtain an expression for
$d_{\alpha}^{\ell}(\theta,y^{(n)}) + d_{\alpha}^u(\theta,y^{(n)})$ in terms of $b(\theta, y^{(n)})_i$ 
and $i^{rs}$.
Now, the length of
$I_{\alpha}^+(\widehat \Theta, y^{(n)})$
is equal to
\begin{align*}
&u_{\alpha}(\widehat \Theta,y^{(n)}) - \ell_{\alpha}(\widehat \Theta,y^{(n)}) + d_{\alpha}^u(\widehat \Theta,y^{(n)})
 + d_{\alpha}^{\ell}(\widehat \Theta,y^{(n)}) \\
&=u_{\alpha} (\theta,y^{(n)}) + 
\frac{\partial u_{\alpha} (\theta,y^{(n)})}{\partial \theta_i} (\widehat \Theta_i - \theta_i)
+ \frac{1}{2} \frac{\partial^2 u_{\alpha} (\theta,y^{(n)})}{\partial \theta_r \partial \theta_s}
(\widehat \Theta_r - \theta_r)(\widehat \Theta_s - \theta_s) \\
& \phantom{123}- \ell_{\alpha} (\theta,y^{(n)}) - 
\frac{\partial \ell_{\alpha} (\theta,y^{(n)})}{\partial \theta_i} (\widehat \Theta_i - \theta_i)
- \frac{1}{2} \frac{\partial^2 \ell_{\alpha} (\theta,y^{(n)})}{\partial \theta_r \partial \theta_s}
(\widehat \Theta_r - \theta_r)(\widehat \Theta_s - \theta_s) \\
& \phantom{123} + d_{\alpha}^u(\theta,y^{(n)}) + d_{\alpha}^{\ell}(\theta,y^{(n)}) + \cdots 
\end{align*}
Thus $E_{\theta} \big(\text{length of } I_{\alpha}^+(\widehat \Theta, y^{(n)}) \, \big| \, Y^{(n)} = y^{(n)}\big)$ 
is equal to
\begin{align*}
u_{\alpha} (\theta,y^{(n)}) - \ell_{\alpha} (\theta,y^{(n)}) 
- \frac{1}{2} &\Bigg( \frac{f^{\prime}(u_{\alpha} (\theta,y^{(n)}); \theta,y^{(n)})}
{f(u_{\alpha} (\theta,y^{(n)}); \theta,y^{(n)})} 
\frac{\partial u_{\alpha} (\theta,y^{(n)})}{\partial \theta_r} 
\frac{\partial u_{\alpha} (\theta,y^{(n)})}{\partial \theta_s} \\
&- \frac{f^{\prime}(\ell_{\alpha} (\theta,y^{(n)}); \theta,y^{(n)})}
{f(\ell_{\alpha} (\theta,y^{(n)}); \theta,y^{(n)})} 
\frac{\partial \ell_{\alpha} (\theta,y^{(n)})}{\partial \theta_r} 
\frac{\partial \ell_{\alpha} (\theta,y^{(n)})}{\partial \theta_s} \Bigg )
i^{rs} + \cdots
\end{align*}
We see that the asymptotic conditional bias $b(\theta, y^{(n)})n^{-1}$ does not enter into this 
expression. 
This has the following two consequences. Firstly, if the improved prediction interval is 
found algebraically using \eqref{c} and \eqref{sum_2_d} then we can use that estimator $\widehat \Theta$ whose 
asymptotic conditional bias is easiest to find. Usually, this will be the (conditional) maximum likelihood estimator
whose asymptotic conditional bias can be found using the formula of
 Vidoni (2004, p.144). Secondly, if
the simulation-based method, similar to that of Kabaila and Syuhada (2008),
is used then we know that the asymptotic efficiency of 
the improved $1-\alpha$ prediction limit is 
independent of the estimator $\widehat \Theta$, on which the estimative $1-\alpha$ prediction limit
is based. Note that we assume throughout this paper that $\widehat \Theta$ has the same asymptotic distribution as the 
(conditional) maximum likelihood estimator.

Now consider the particular case that $f(\, \cdot \, ; \theta, y^{(n)})$
is also symmetric about $m(\theta,y^{(n)})$  for all $y^{(n)}$ and $\theta$. In other words,
suppose that, for every $y^{(n)}$ and $\theta$, 
$f \big(m(\theta,y^{(n)}) - w ; \theta, y^{(n)} \big)=
f \big(m(\theta,y^{(n)}) + w ; \theta, y^{(n)}\big)$ for all $w > 0$. In this case, 
we may choose 
$d_{\alpha}^{\ell}(\theta,y^{(n)}) = d_{\alpha}^u(\theta,y^{(n)}) = \delta_{\alpha}(\theta,y^{(n)})$,
say.
Define $w_{\alpha}(\theta,y^{(n)})$  by the requirement that
$P_\theta \big ( Y_{n+k} \in \big [ m(\theta,y^{(n)}) - w_{\alpha}(\theta,y^{(n)}), 
\, m(\theta,y^{(n)}) + w_{\alpha}(\theta,y^{(n)}) \big ] \, \big| \,
Y^{(n)} = y^{(n)} \big ) = 1 - \alpha$
for all $\theta$ and $y^{(n)}$. Thus 
$\ell_{\alpha}(\theta,y^{(n)}) =  m_{\alpha}(\theta,y^{(n)}) - w_{\alpha}(\theta,y^{(n)})$
and 
$u_{\alpha}(\theta,y^{(n)}) =  m_{\alpha}(\theta,y^{(n)}) + w_{\alpha}(\theta,y^{(n)})$.
It may be shown that $\delta_{\alpha}(\theta,y^{(n)})$ is equal to
\begin{align*}
- n^{-1} \, \frac{\partial w_{\alpha}(\theta,y^{(n)})}{\partial \theta_i} \,
b(\theta, y^{(n)})_i  
- \frac{1}{4 f(u_{\alpha}(\theta,y^{(n)});\theta,y^{(n)})} 
\frac{\partial^2 G_{\alpha}(\widehat \theta; \theta|y^{(n)})}{\partial \widehat \theta_r \partial \widehat \theta_s}
\Bigg|_{\hat{\theta} = \theta} i^{rs}
\end{align*}
The improved prediction interval $\big [ \ell_{\alpha}(\widehat \Theta,Y^{(n)}) - \delta_{\alpha}(\widehat \Theta,Y^{(n)}), 
\, u_{\alpha}(\widehat \Theta,Y^{(n)}) + \delta_{\alpha}(\widehat \Theta,Y^{(n)}) \big ]$
has been considered in the context of one-step-ahead prediction for $\{Y_t\}$ a 
stationary zero-mean Gaussian AR(1) process by Kabaila and Syuhada (2007, Section 4), where the formula
for $d(\theta, y_n)$ should be $-c(\theta, y_n)/(2 v^{-1/2} \phi(z_{1-\frac{\alpha}{2}}))$ instead of
$-c(\theta, y_n)/(2 \phi(z_{1-\frac{\alpha}{2}}))$.


\bigskip
\begin{center}
\large{REFERENCES}
\end{center}

\noindent BARNDORFF-NIELSEN, O.E. and COX, D.R. (1994) {\it
Inference and Asymptotics.} London: Chapman and Hall.

\noindent CHRISTOFFERSEN, P.F. (1998) Evaluating interval forecasts.
{\it International Economic Review}  39, 841-862.

\noindent HE, Z. (2000) Assessment of the accuracy of time series predictions.
Unpublished Ph.D. thesis, Department of Statistical Science, La Trobe University.

\noindent KABAILA, P. (1993) On bootstrap predictive inference for autoregressive processes.
{\it Journal of Time Series Analysis} 14, 473-484.

\noindent KABAILA, P. (1999) The relevance property for prediction intervals.
{\it Journal of Time Series Analysis} 20, 655-662.

\noindent KABAILA, P. and HE, Z. (2004) The adjustment of
prediction intervals to account for errors in parameter
estimation. {\it Journal of Time Series Analysis} 25, 351-358.

\noindent KABAILA, P. and SYUHADA, K. (2007) The relative
efficiency of prediction intervals. {\it Communications in Statistics -
Theory and Methods} 36, 2673--2686.

\noindent KABAILA, P. and SYUHADA, K. (2008) Improved prediction limits for AR($p$) and
ARCH($p$) processes. {\it Journal of Time Series Analysis} 29, 213--223.

\noindent McCULLOUGH, B.D. (1994) Bootstrapping forecast intervals: an application to AR(p) models.
{\it Journal of Forecasting} 13, 51-66.

\noindent PHILLIPS, P.C.B. (1979) The sampling distribution of forecasts from a first-order
autoregression. {\it Journal of Econometrics} 9, 241-261.

\noindent SHAMAN, P. and STINE, R.A. (1988) The bias of
autoregressive coefficient estimators. {\it Journal of the
American Statistical Association} 83, 842-848.

\noindent STINE, R.A. (1987) Estimating properties of autoregressive forecasts.
{\it Journal of the American Statistical Association} 82, 1072-1078.

\noindent THOMBS, L.A. and SCHUCANY, W.R. (1990) Bootstrap prediction intervals for
autoregression.
{\it Journal of the American Statistical Association} 85, 486-492.

\noindent VIDONI, P. (2004) Improved prediction intervals for
stochastic process models. {\it Journal of Time Series Analysis}
25, 137-154

\end{document}